\documentclass[12pt]{article}
\usepackage{latexsym, amssymb, amsmath, amscd, amsfonts, epsfig, graphicx, colordvi,verbatim,ifpdf,extarrows}
\usepackage{amsfonts, amsmath, amssymb}
\usepackage{amssymb,amsfonts,amsmath,latexsym,epsfig,cite, psfrag,eepic,color}
\usepackage{amscd,graphics}
\usepackage{latexsym, amssymb,  amsmath,amscd, amsfonts, epsfig, graphicx, colordvi,amsthm}

\usepackage{graphicx}
\usepackage{epstopdf}
\usepackage{color}
\usepackage{ifpdf}
\usepackage{fancybox}
\usepackage[font=small,labelfont=bf,labelsep=none]{caption}
\usepackage{float}

\allowdisplaybreaks

\newtheorem{thm}{Theorem}[section]

\newtheorem{defi}[thm]{Definition}
\newtheorem{lem}[thm]{Lemma}

\def\pf{\noindent{\it Proof.} }
\def\qed{\nopagebreak\hfill{\rule{4pt}{7pt}}
\medbreak}

\numberwithin{equation}{section}

\def\qed{\nopagebreak\hfill{\rule{4pt}{7pt}}
\medbreak}

\newlength{\boxedparwidth}
  {\begin{center} \begin{tabular}{|@{\hspace{.315in}}c@{\hspace{.15in}}|}
                  \hline \\ \begin{minipage}[t]{\boxedparwidth}
                  \setlength{\parindent}{.25in}}%
  {\end{minipage} \\ \\ \hline \end{tabular} \end{center}}

\begin{document}

\begin{center}
{\Large \bf Six types of separable integer partitions}
\end{center}

\begin{center}
 {Thomas Y. He}$^{1}$, {Y. Hu}$^{2}$,  {H.X. Huang}$^{3}$ and {Y.X. Xie}$^{4}$ \vskip 2mm

$^{1,2,3,4}$ School of Mathematical Sciences, Sichuan Normal University, Chengdu 610066, P.R. China

   \vskip 2mm

 $^1$heyao@sicnu.edu.cn,  $^2$huyue@stu.sicnu.edu.cn, $^3$huanghaoxuan@stu.sicnu.edu.cn, $^4$xieyx@stu.sicnu.edu.cn
\end{center}

\vskip 6mm   {\noindent \bf Abstract.} Recently, Andrews introduced separable integer partition classes and studied some well-known theorems. In this atricle, we will investigate six types of partitions from the view of the point of  separable integer partition classes.

\noindent {\bf Keywords}: separable integer partition classes, the largest part, generating functions

\noindent {\bf AMS Classifications}: 05A17, 11P83

\section{Introduction}

A partition $\pi$ of a positive integer $n$ is a finite non-increasing sequence of positive integers $\pi=(\pi_1,\pi_2,\ldots,\pi_\ell)$ such that $\pi_1+\pi_2+\cdots+\pi_\ell=n$. The empty sequence forms the only partition of zero. The $\pi_i$ are called the parts of $\pi$. Let $\ell(\pi)$ be the number of parts of $\pi$ and let $l(\pi)$ be the largest part of $\pi$. The weight of $\pi$ is the sum of parts, denoted $|\pi|$.

In \cite{Andrews-2022}, Andrews introduced separable integer partition classes with modulus $k$ and  studied some well-known theorems from the point of view of separable integer partition classes with modulus $k$.

 \begin{defi}
For a positive integer $k$, a separable integer partition class $\mathcal{P}$ with modulus $k$ is a subset of all the partitions satisfying the following{\rm:}

There is a subset $\mathcal{B}\subset\mathcal{P}$ {\rm(}$\mathcal{B}$ is called the basis of $\mathcal{P}${\rm)} such that for each integer $m\geq 1$, the number of partitions in $\mathcal{B}$ with $m$ parts is finite and every partition in $\mathcal{P}$ with $m$ parts is uniquely of the form
\begin{equation}\label{ordinary-form-1}
(b_1+\pi_1)+(b_2+\pi_2)+\cdots+(b_m+\pi_m),
\end{equation}
where $(b_1,b_2,\ldots,b_m)$ is a partition in $\mathcal{B}$ and $(\pi_1,\pi_2,\ldots,\pi_m)$ is a  non-increasing sequence of nonnegative integers, whose only restriction is that each part is divisible by $k$. Furthermore, all partitions of the form \eqref{ordinary-form-1} are in $\mathcal{P}$.
\end{defi}

Some separable integer partition classes are studied in \cite{Chen-He-Huang-Zhang-2025,Chen-He-Tang-Wei-2024,He-Huang-Li-Zhang-2025,Passary-2019}. In this article, we will investigate six types of partitions with the aid of  separable integer partition classes with modulus $k$.

For integers $k,r$ and $N$ such that $k\geq r\geq 1$ and $N\geq 1$, we define $\varphi_{k,r}(N)$ to be the integer $s$ such that $1\leq s\leq r$ and $N\equiv s\pmod{k}$. Now, we introduce two sets of partitions, $\mathcal{P}_{k,r}$ and $\mathcal{D}_{k,r}$.

\begin{defi}
For integers $k$ and $r$ such that $k\geq r\geq 1$, let $\mathcal{P}_{k,r}$ denote the set of partitions $(\pi_1,\pi_2,\ldots,\pi_\ell)$ such that
\begin{itemize}
\item[\rm(1)] for $1\leq i\leq \ell$, $1\leq\varphi_{k,r}(\pi_i)\leq r$;

\item[\rm(2)] for $1\leq i<\ell$, $\varphi_{k,r}(\pi_{i})\geq\varphi_{k,r}(\pi_{i+1})$.
\end{itemize}

\end{defi}

\begin{defi}
For integers $k$ and $r$ such that $k\geq r\geq 1$, let $\mathcal{D}_{k,r}$ denote the set of partitions $(\pi_1,\pi_2,\ldots,\pi_\ell)$ such that
\begin{itemize}
\item[\rm(1)] for $1\leq i\leq \ell$, $1\leq\varphi_{k,r}(\pi_i)\leq r$;

\item[\rm(2)] for $1\leq i<\ell$, $\varphi_{k,r}(\pi_{i})\geq\varphi_{k,r}(\pi_{i+1})$ and  $\pi_{i}-\pi_{i+1}\geq k$.
\end{itemize}

\end{defi}

If $k=r=1$, then  $\mathcal{P}_{1,1}$ is the set of all partitions and $\mathcal{D}_{1,1}$ is the set of all distinct partitions. For $k\geq r\geq 1$, we will give the generating functions for the partitions in $\mathcal{P}_{k,r}$ and $\mathcal{D}_{k,r}$ respectively with the aid of separable integer partition classes with modulus $k$.
\begin{equation}\label{gen-p-kr}
\sum_{\pi\in\mathcal{P}_{k,r}}z^{\ell(\pi)}q^{|\pi|}=\sum_{m\geq 0}\frac{z^mq^m}{(q^k;q^k)_m}{{m+r-1}\brack{r-1}}_1,
\end{equation}
and
\begin{equation}\label{gen-d-kr}
\sum_{\pi\in\mathcal{D}_{k,r}}z^{\ell(\pi)}q^{|\pi|}=\sum_{m\geq 0}\frac{z^mq^{k{m\choose 2}+m}}{(q^k;q^k)_m}{{m+r-1}\brack{r-1}}_1.
\end{equation}
Here and in the sequel, we assume that $|q|<1$ and employ the standard notation \cite{Andrews-1976}:
\[(a;q)_\infty=\prod_{i\geq 0}(1-aq^i),\]
\[(a;q)_n=\frac{(a;q)_\infty}{(aq^n;q)_\infty},\]
and ${A\brack B}_k$ is the $q$-binomial coefficient, or Gaussian polynomial for nonnegative integers $A$ and $B$ defined as follows:
\[{A\brack B}_k=\left\{\begin{array}{ll}\frac{(q^k;q^k)_A}{(q^k;q^k)_B(q^k;q^k)_{A-B}},&\text{if }A\geq B\geq 0,\\
0,&\text{otherwise.}
\end{array}
\right.\]

Assume that $a,b$ and $k$ are integers such that $k\geq b>a\geq1$. Throughout this article, we use $\ell_{a}(\pi)$ and $\ell_{b}(\pi)$ to denote the number of parts equivalent to $a$ and $b$ modulo $k$ in a partition $\pi$ respectively. Let $\mathcal{P}_{a,b,k}$ be the set of partitions whose parts are equivalent to $a$ or $b$ modulo $k$. Clearly, the generating function for the partitions in $\mathcal{P}_{a,b,k}$ is
 \[\sum_{\pi\in\mathcal{P}_{a,b,k}}\mu^{\ell_{a}(\pi)}\nu^{\ell_{b}(\pi)}q^{|\pi|}=\frac{1}{(\mu q^a;q^k)_\infty(\nu q^b;q^k)_\infty}.\]
 In  \cite{He-Huang-Li-Zhang-2025}, He, Huang, Li and Zhang showed that $\mathcal{P}_{a,b,k}$ is a separable integer partition class with modulus $k$ and gave a new generating function for the partitions in $\mathcal{P}_{a,b,k}$.
\begin{align*}
&\sum_{\pi\in\mathcal{P}_{a,b,k}}\mu^{\ell_{a}(\pi)}\nu^{\ell_{b}(\pi)}q^{|\pi|}\\
&=\sum_{m,h,i\geq 0}\frac{\mu^{m-h-i}\nu^{h+i}q^{ma+kh^2+(b-a)(h+i)}}{(q^k;q^k)_m}{{h+i}\brack{h}}_k{{m-h-i}\brack{h}}_k.
\end{align*}

In this article, we will consider four subsets of $\mathcal{P}_{a,b,k}$.
We first introduce two subsets $\mathcal{R}_{a,b,k}$ and $\mathcal{R}'_{a,b,k}$ of $\mathcal{P}_{a,b,k}$.
\begin{defi}
Assume that $a,b$ and $k$ are integers such that $k\geq b>a\geq1$. Let $\mathcal{R}_{a,b,k}$ (resp. $\mathcal{R}'_{a,b,k}$) denote the set of partitions $(\pi_1,\pi_2,\ldots,\pi_\ell)$ such that
\begin{itemize}
\item[(1)] for $1\leq i\leq \ell$, $\pi_i\equiv a$ or $b\pmod{k}$;

\item[(2)] only parts equivalent to $a$ (resp. $b$) modulo $k$ may be repeated.
\end{itemize}
\end{defi}

Clearly, we have
\begin{equation}\label{gen-r-abk}
\sum_{\pi\in\mathcal{R}_{a,b,k}}u^{\ell_a(\pi)}v^{\ell_b(\pi)}q^{|\pi|}=\frac{(-vq^{b};q^k)_\infty}{(uq^a;q^k)_\infty},
\end{equation}
and
\begin{equation}\label{gen-rr-abk}
\sum_{\pi\in\mathcal{R}'_{a,b,k}}u^{\ell_a(\pi)}v^{\ell_b(\pi)}q^{|\pi|}=\frac{(-uq^{a};q^k)_\infty}{(vq^b;q^k)_\infty}.
\end{equation}

We will show that $\mathcal{R}_{a,b,k}$ and $\mathcal{R}'_{a,b,k}$ are separable integer partition classes with modulus $k$ and then give the proofs of \eqref{gen-r-abk} and \eqref{gen-rr-abk}.

Then, we introduce two subsets $\mathcal{D}_{a,b,k}$ and $\mathcal{D}'_{a,b,k}$ of $\mathcal{P}_{a,b,k}$.
\begin{defi}
Assume that $a,b$ and $k$ are integers such that $k\geq b>a\geq1$. Let $\mathcal{D}_{a,b,k}$ (resp. $\mathcal{D}'_{a,b,k}$) denote the set of partitions $(\pi_1,\pi_2,\ldots,\pi_\ell)$ such that
\begin{itemize}
\item[(1)] for $1\leq i\leq \ell$, $\pi_i\equiv a$ or $b\pmod{k}$;

\item[(2)] for $1\leq i<\ell$, $\pi_i-\pi_{i+1}\geq k$ with strict inequality if $\pi_i\equiv b\pmod{k}$ (resp. $\pi_i\equiv a\pmod{k}$).
\end{itemize}
\end{defi}

If $a=1,b=2$ and $k=2$, then $\mathcal{D}_{1,2,2}$ is the set of partitions whose parts satisfy the difference conditions in the first G\"ollnitz-Gordon identity \cite{Gollnitz-1967,Gordon-1965} and $\mathcal{D}'_{1,2,2}$ is the set of partitions whose parts satisfy the difference conditions in the first little G\"ollnitz  identity \cite{Gollnitz-1967}. The first G\"ollnitz-Gordon identity and the first little G\"ollnitz  identity have been studied in view of separable integer partition classes with modulus $k$ by Andrews \cite{Andrews-2022} and Passary \cite{Passary-2019} respectively.

For $k\geq b>a\geq1$,
an easy combinatorial argument implies that
\begin{equation}\label{gen-d-abk}
\begin{split}
\sum_{\pi\in\mathcal{D}_{a,b,k}}u^{\ell_a(\pi)}v^{\ell_b(\pi)}q^{|\pi|}&=\sum_{m\geq 0}\frac{u^mq^{k{m\choose 2}+ma}}{(q^k;q^k)_m}(-u^{-1}vq^{b-a};q^k)_m\\
&=\sum_{h\geq0}\frac{v^hq^{2k{h\choose 2}+bh}}{(q^k;q^k)_h}(-uq^{kh+a};q^k)_\infty,
\end{split}
\end{equation}
and
\begin{equation}\label{gen-dd-abk}
\begin{split}
\sum_{\pi\in\mathcal{D'}_{a,b,k}}u^{\ell_a(\pi)}v^{\ell_b(\pi)}q^{|\pi|}&=\sum_{m\geq 0}\frac{v^mq^{k{m\choose 2}+mb}}{(q^k;q^k)_m}(-uv^{-1}q^{a-b};q^k)_m\\
&=\sum_{h\geq 0}\frac{u^hq^{2k{h\choose 2}+ah}}{(q^k;q^k)_h}(-vq^{kh+b};q^k)_\infty.
\end{split}
\end{equation}

For $\ell\geq 0$, let $\mathcal{R}_{a,b,k}(\ell)$, $\mathcal{R}'_{a,b,k}(\ell)$, $\mathcal{D}_{a,b,k}(\ell)$ and $\mathcal{D}'_{a,b,k}(\ell)$ be the set of partitions in $\mathcal{R}_{a,b,k}$, $\mathcal{R}'_{a,b,k}$, $\mathcal{D}_{a,b,k}$ and $\mathcal{D}'_{a,b,k}$ with $\ell$ parts respectively.
For a partition $\pi=(\pi_1,\pi_2,\ldots,\pi_{\ell-1},\pi_\ell)$ in $\mathcal{R}_{a,b,k}(\ell)$ (resp. $\mathcal{R}'_{a,b,k}(\ell)$), if we add $k(\ell-1),k(\ell-2),\ldots,k$ to the parts $\pi_1,\pi_2,\ldots,\pi_{\ell-1}$ respectively, then we get a partition
 $\lambda=(\pi_1+k(\ell-1),\pi_2+k(\ell-2),\ldots,\pi_{\ell-1}+k,\pi_\ell)$ in $\mathcal{D}_{a,b,k}(\ell)$ (resp. $\mathcal{D}'_{a,b,k}(\ell)$), and vice versa. Moreover, we have $\ell_a(\pi)=\ell_a(\lambda)$ and  $\ell_b(\pi)=\ell_b(\lambda)$. This implies that
 \[\sum_{\pi\in\mathcal{D}_{a,b,k}(\ell)}u^{\ell_a(\pi)}v^{\ell_b(\pi)}q^{|\pi|}=q^{k{\ell\choose 2}}\sum_{\pi\in\mathcal{R}_{a,b,k}(\ell)}u^{\ell_a(\pi)}v^{\ell_b(\pi)}q^{|\pi|},\]
 and
 \[\sum_{\pi\in\mathcal{D}'_{a,b,k}(\ell)}u^{\ell_a(\pi)}v^{\ell_b(\pi)}q^{|\pi|}=q^{k{\ell\choose 2}}\sum_{\pi\in\mathcal{R}'_{a,b,k}(\ell)}u^{\ell_a(\pi)}v^{\ell_b(\pi)}q^{|\pi|}.\]

For completeness, we will show \eqref{gen-d-abk} and \eqref{gen-dd-abk} from the view of the point of separable integer partition classes with modulus $k$.

This article is organized as follows. In Section 2, we recall some necessary identities. With the aid of separable integer partition classes with modulus $k$, we give the proofs of \eqref{gen-p-kr} and \eqref{gen-d-kr},  \eqref{gen-r-abk} and \eqref{gen-rr-abk}, and \eqref{gen-d-abk} and \eqref{gen-dd-abk} in Section 3, Section 4 and Section 5  respectively.

\section{Preliminaries}

In this section, we collect some classical identities needed in this article from \cite{Andrews-1976}.

{\noindent \bf The $q$-binomial theorem \cite[Theorem 2.1]{Andrews-1976}:}
\begin{equation}\label{chang-1}
\sum_{n\geq0} \frac{(a;q)_{n}}{(q;q)_{n}}t^{n}=\frac{{(at;q)}_{\infty}}{(t;q)_{\infty}}.
\end{equation}

Euler found the following two special cases of \eqref{chang-1} (see \cite[Corollary 2.2]{Andrews-1976}):
\begin{equation}\label{Euler-1}
\sum_{m\geq0}\frac{t^m}{(q;q)_m}=\frac{1}{(t;q)_\infty},
\end{equation}
and
\begin{equation}\label{Euler-2}
\sum_{m\geq0}\frac{t^mq^{{m}\choose 2}}{(q;q)_m}=(-t;q)_\infty.
\end{equation}

The following two recurrences for the $q$-binomial coefficients are needed.

{\noindent\bf \cite[(3.3.4)]{Andrews-1976}:}
\begin{equation}\label{bin-new-r-1}
{A\brack B}_k={{A-1}\brack{B-1}}_k+q^{kB}{{A-1}\brack{B}}_k\text{ for }A\geq B\geq 1.
\end{equation}

{\noindent\bf \cite[(3.3.9)]{Andrews-1976}:}
\begin{equation}\label{bin-new-r-s}{{A+B+1}\brack{B+1}}_1=\sum_{s=0}^{A}q^s{{B+s}\brack{B}}_1\text{ for }A,B\geq 0.
\end{equation}

We need the following formula related to the $q$-binomial coefficients \cite[(3.3.6)]{Andrews-1976}.
\begin{equation}\label{h-sum}
(t;q)_{m}=\sum_{h=0}^{m}{m\brack h}(-t)^{h}q^{\binom{h}{2}}.
\end{equation}

\section{Proofs of \eqref{gen-p-kr} and \eqref{gen-d-kr}}

In this section, we assume that $k$ and $r$ are integers such that $k\geq r\geq 1$. Let $\mathcal{BP}_{k,r}$ be the set of partitions  in $\mathcal{P}_{k,r}$ with parts not exceeding $r$ and let $\mathcal{BD}_{k,r}$ be the set of partitions $(\lambda_1,\lambda_2,\ldots,\lambda_\ell)$ in $\mathcal{D}_{k,r}$ such that
\begin{itemize}
\item[(1)] $1\leq\lambda_\ell\leq r$;

\item[(2)] for $1\leq i<\ell$, $\lambda_i-\lambda_{i+1}<2k$.
\end{itemize}

Clearly, $\mathcal{BP}_{k,r}$ and $\mathcal{BD}_{k,r}$ are the basis of $\mathcal{P}_{k,r}$ and  $\mathcal{D}_{k,r}$ respectively. So we have

\begin{lem}
$\mathcal{P}_{k,r}$ and  $\mathcal{D}_{k,r}$ are separable integer partition classes with modulus $k$.
\end{lem}

For $m\geq1$, let $\mathcal{BP}_{k,r}(m)$ (resp. $\mathcal{BD}_{k,r}(m)$) be the set of partitions in $\mathcal{BP}_{k,r}$ (resp. $\mathcal{BD}_{k,r}$) with $m$ parts. For a partition $\lambda=(\lambda_1,\lambda_2,\ldots,\lambda_{m-1},\lambda_m)$ in $\mathcal{BP}_{k,r}(m)$, if we add $k(m-1),k(m-2),\ldots,k$ to the parts $\lambda_1,\lambda_2,\ldots,\lambda_{m-1}$ respectively, then we get a partition
 $(\lambda_1+k(m-1),\lambda_2+k(m-2),\ldots,\lambda_{m-1}+k,\lambda_m)$ in $\mathcal{BD}_{k,r}(m)$, and vice versa. This implies that
 \[\sum_{\lambda\in\mathcal{BD}_{k,r}(m)}q^{|\lambda|}=q^{k{m\choose 2}}\sum_{\lambda\in\mathcal{BP}_{k,r}(m)}q^{|\lambda|}.\]

We find that in order to prove \eqref{gen-p-kr} and \eqref{gen-d-kr}, it suffices to show the following theorem.
\begin{thm}For $m\geq 1$,
\[\sum_{\lambda\in\mathcal{BP}_{k,r}(m)}q^{|\lambda|}=q^m{{m+r-1}\brack{r-1}}_1.\]
\end{thm}

\pf For $m\geq1$ and $1\leq s\leq r$, let $\mathcal{BP}_{k,r}(m,s)$ be the set of partitions in $\mathcal{BP}_{k,r}(m)$ with the largest part being $s$. Assume that $\lambda$ is a partition in $\mathcal{BP}_{k,r}(m,s)$, that is, $\ell(\lambda)=m$ and $l(\lambda)=s$.
If we remove the largest part $s$ from $\lambda$ and subtract $1$ from each of the remaining parts, then we get a partition with at most $m-1$ parts which do not exceed $s-1$, and vice versa. It yields that
\[\sum_{\lambda\in\mathcal{BP}_{k,r}(m,s)}q^{|\lambda|}=q^{m+s-1}{{m+s-2}\brack{s-1}}_1.\]

Then, we obtain that for $m\geq 1$,
\begin{align}
\sum_{\lambda\in\mathcal{BP}_{k,r}(m)}q^{|\lambda|}
&=\sum_{s=1}^r\sum_{\lambda\in\mathcal{BP}_{k,r}(m,s)}q^{|\lambda|}\nonumber\\
&=\sum_{s=1}^rq^{m+s-1}{{m+s-2}\brack{s-1}}_1\nonumber\\
&=q^m\sum_{s=0}^{r-1}q^{s}{{m-1+s}\brack{s}}_1\nonumber\\
&=q^m{{m+r-1}\brack{r-1}}_1,\nonumber
\end{align}
where the final equation follows from \eqref{bin-new-r-s} with $A=r-1$ and $B=m-1$. The proof is complete. \qed

\section{Proofs of \eqref{gen-r-abk} and \eqref{gen-rr-abk}}

In the remainder of this article, we assume that $a,b$ and $k$ are integers such that $k\geq b>a\geq1$.
In this section, we will show \eqref{gen-r-abk} and \eqref{gen-rr-abk} in consideration of separable integer partition classes with modulus $k$ in Section 4.1 and Section 4.2 respectively.

\subsection{Proof of \eqref{gen-r-abk}}

We first  show that $\mathcal{R}_{a,b,k}$ is a separable integer partition class with modulus $k$. To do this, we are required to find the basis of $\mathcal{R}_{a,b,k}$, which involves the following set. For $m\geq 1$, let $\mathcal{BR}_{a,b,k}(m)$ be the set of partitions $\lambda=(\lambda_1,\lambda_2,\ldots,\lambda_m)$ in $\mathcal{R}_{a,b,k}$ with $m$ parts
 such that
\begin{itemize}
\item[(1)] $\lambda_m=a$ or $b$;

\item[(2)] for $1\leq i<m$, $\lambda_i-\lambda_{i+1}\leq k$ with strict inequality if $\lambda_{i+1}\equiv a\pmod{k}$.
\end{itemize}

Assume that $\lambda$ is a partition in $\mathcal{BR}_{a,b,k}(m)$. For $1\leq i<m$, if $\lambda_{i+1}=kj+a$, then we have $kj+a\leq\lambda_{i}<k(j+1)+a$, and so  $\lambda_{i}=kj+a$ or $kj+b$. If $\lambda_{i+1}=kj+b$, then we have $kj+b<\lambda_{i}\leq k(j+1)+b$, and so  $\lambda_{i}=k(j+1)+a$ or $k(j+1)+b$. So, the number of partitions in $\mathcal{BR}_{a,b,k}(m)$ is $2^m$. For example, the number of partitions in $\mathcal{BR}_{a,b,k}(3)$ is $2^3=8$.
\[(a,a,a),(b,a,a),(k+a,b,a),(k+b,b,a),\]
\[(k+a,k+a,b),(k+b,k+a,b),(2k+a,k+b,b),(2k+b,k+b,b).\]

\begin{lem}
 $\mathcal{R}_{a,b,k}$ is a separable integer partition class with modulus $k$.
\end{lem}

\pf  Set \[\mathcal{BR}_{a,b,k}=\bigcup_{m\geq 1}\mathcal{BR}_{a,b,k}(m).\]
 Obviously, $\mathcal{BR}_{a,b,k}$ is the basis of $\mathcal{R}_{a,b,k}$. This completes the proof.   \qed

For $m\geq 1$, let $\lambda$ be a partition in $\mathcal{BR}_{a,b,k}(m)$. Clearly, we have
\[a\leq l(\lambda)\leq k(m-1)+b.\]

For $m\geq 1$ and $0\leq h\leq m-1$, let $\mathcal{BR}_{a,b,k}(m,h,a)$ (resp. $\mathcal{BR}_{a,b,k}(m,h,b)$) be the set of partitions in $\mathcal{BR}_{a,b,k}(m)$ with the largest part $kh+a$ (resp. $kh+b$). The generating functions for the partitions in $\mathcal{BR}_{a,b,k}(m,h,a)$ and $\mathcal{BR}_{a,b,k}(m,h,b)$ respectively are given below.
\begin{thm}\label{thm-gen-r-abk-basis}
For $m\geq 1$ and $0\leq h\leq m-1$, we have
\begin{equation*}\label{eqn-gen-r-abk-basis-a}
\sum_{\lambda\in\mathcal{BR}_{a,b,k}(m,h,a)}u^{\ell_a(\lambda)}v^{\ell_b(\lambda)}q^{|\lambda|}=u^{m-h}v^{h}q^{ma+k{{h+1}\choose 2}+(b-a)h}{{m-1}\brack{h}}_k,
\end{equation*}
and
\begin{equation}\label{eqn-gen-rr-abk-basis-b}
\sum_{\lambda\in\mathcal{BR}_{a,b,k}(m,h,b)}u^{\ell_a(\lambda)}v^{\ell_b(\lambda)}q^{|\lambda|}
=u^{m-h-1}v^{h+1}q^{ma+k{{h+1}\choose 2}+(b-a)(h+1)}{{m-1}\brack{h}}_k.
\end{equation}
\end{thm}

\pf Let $\lambda$ be a partition in $\mathcal{BR}_{a,b,k}(m,h,b)$. By definition, we have $l(\lambda)=kh+b$. If we change the largest part $kh+b$ of $\lambda$ to $kh+a$, then we get a partition in $\mathcal{BR}_{a,b,k}(m,h,a)$, and vice versa. This implies that
\[\sum_{\lambda\in\mathcal{BR}_{a,b,k}(m,h,b)}u^{\ell_a(\lambda)}v^{\ell_b(\lambda)}q^{|\lambda|}=u^{-1}vq^{b-a}\sum_{\lambda\in\mathcal{BR}_{a,b,k}(m,h,a)}u^{\ell_a(\lambda)}v^{\ell_b(\lambda)}q^{|\lambda|}.\]  We find that it suffices to show \eqref{eqn-gen-rr-abk-basis-b}.

For a partition $\lambda$ in $\mathcal{BR}_{a,b,k}(m,h,b)$, it is clear that there exist parts
\[kh+b,k(h-1)+b,\ldots,b\]
in $\lambda$. If we remove the parts $kh+b,k(h-1)+b,\ldots,b$ from $\lambda$, then we get a partition
with $m-h-1$ parts not exceeding ${kh+a}$ and equivalent to $a$ modulo $k$. The generating function for the partitions with $m-h-1$ parts not exceeding ${kh+a}$ and equivalent to $a$ modulo $k$ is
\[q^{a(m-h-1)}{{m-1}\brack{h}}_k.\]

So, we get
\begin{align*}
\sum_{\lambda\in\mathcal{BR}_{a,b,k}(m,h,b)}u^{\ell_a(\lambda)}v^{\ell_b(\lambda)}q^{|\lambda|}=v^{h+1}q^{k{{h+1}\choose 2}+b(h+1)}\cdot u^{m-h-1}q^{a(m-h-1)}{{m-1}\brack{h}}_k.
\end{align*}
Hence,  \eqref{eqn-gen-rr-abk-basis-b} is verified. The proof is complete.    \qed

Now, we proceed to give a proof of \eqref{gen-r-abk}.

{\noindent \bf Proof of \eqref{gen-r-abk}.} We first show that for $m\geq 1$,
\begin{equation}\label{000000000000001}
\sum_{\lambda\in\mathcal{BR}_{a,b,k}(m)}u^{\ell_a(\lambda)}v^{\ell_b(\lambda)}q^{|\lambda|}=\sum_{h=0}^{m}u^{m-h}v^{h}q^{ma+k{h\choose 2}+(b-a)h}{{m}\brack{h}}_k.
\end{equation}

 For $m=1$, since there are two partitions $(a)$ and $(b)$ in $\mathcal{BR}_{a,b,k}(1)$, we have
\begin{equation*}
\sum_{\lambda\in\mathcal{BD}'_{a,b,k}(1)}u^{\ell_a(\lambda)}v^{\ell_b(\lambda)}q^{|\lambda|}=uq^a+vq^b,
\end{equation*}
which agrees with \eqref{000000000000001} for $m=1$.

For $m\geq 2$, by Theorem \ref{thm-gen-r-abk-basis}, we get
\begin{align*}
&\quad\sum_{\lambda\in\mathcal{BR}_{a,b,k}(m)}u^{\ell_a(\lambda)}v^{\ell_b(\lambda)}q^{|\lambda|}\\
&=\sum_{h=0}^{m-1}\left(\sum_{\lambda\in\mathcal{BR}_{a,b,k}(m,h,a)}u^{\ell_a(\lambda)}v^{\ell_b(\lambda)}q^{|\lambda|}
+\sum_{\lambda\in\mathcal{BR}_{a,b,k}(m,h,b)}u^{\ell_a(\lambda)}v^{\ell_b(\lambda)}q^{|\lambda|}\right)\\
&=\sum_{h=0}^{m-1}\left(u^{m-h}v^{h}q^{ma+k{{h+1}\choose 2}+(b-a)h}{{m-1}\brack{h}}_k+u^{m-h-1}v^{h+1}q^{ma+k{{h+1}\choose 2}+(b-a)(h+1)}{{m-1}\brack{h}}_k\right)\\
&=\sum_{h=0}^{m-1}u^{m-h}v^{h}q^{ma+k{{h+1}\choose 2}+(b-a)h}{{m-1}\brack{h}}_k+\sum_{h=1}^{m}u^{m-h}v^hq^{ma+k{h\choose 2}+(b-a)h}{{m-1}\brack{h-1}}_k\\
&=u^mq^{ma}+\sum_{h=1}^{m-1}u^{m-h}v^{h}q^{ma+k{h\choose 2}+(b-a)h}\left(q^{kh}{{m-1}\brack{h}}_k+{{m-1}\brack{h-1}}_k\right)+v^mq^{ma+k{m\choose 2}+(b-a)m}\\
&=\sum_{h=0}^{m}u^{m-h}v^{h}q^{ma+k{h\choose 2}+(b-a)h}{{m}\brack{h}}_k,
\end{align*}
where the final equation follows from \eqref{bin-new-r-1}. So, we conclude that \eqref{000000000000001} holds for $m\geq 1$.

By virtue of \eqref{000000000000001}, we derive that
\begin{align}
\sum_{\pi\in\mathcal{R}_{a,b,k}}u^{\ell_a(\pi)}v^{\ell_b(\pi)}q^{|\pi|}&=1+\sum_{m\geq 1}\frac{1}{(q^k;q^k)_m}\sum_{\lambda\in\mathcal{BR}_{a,b,k}(m)}u^{\ell_a(\lambda)}v^{\ell_b(\lambda)}q^{|\lambda|}\nonumber\\
&=\sum_{m\geq 0}\frac{1}{(q^k;q^k)_m}\sum_{h=0}^mu^{m-h}v^hq^{ma+k{h\choose 2}+(b-a)h}{m\brack h}_k.\label{0101-1}
\end{align}

Summing up  the $h$-sum in \eqref{0101-1} by letting $q\rightarrow q^{k}$ and $t=-u^{-1}vq^{b-a}$ in \eqref{h-sum}, we get
\begin{align*}
\sum_{\pi\in\mathcal{R}_{a,b,k}}u^{\ell_a(\pi)}v^{\ell_b(\pi)}q^{|\pi|}&=\sum_{m\geq 0}\frac{u^mq^{ma}}{(q^k;q^k)_m}\sum_{h=0}^m(u^{-1}vq^{b-a})^hq^{k{h\choose 2}}{m\brack h}_k\\
&=\sum_{m\geq 0}\frac{u^mq^{ma}}{(q^k;q^k)_m}(-u^{-1}vq^{b-a};q^k)_m\\
&=\frac{(-vq^{b};q^k)_\infty}{(uq^a;q^k)_\infty},
\end{align*}
where the final equation follows from \eqref{chang-1}. We arrive at \eqref{gen-r-abk}.

With the aid of \eqref{Euler-1} and \eqref{Euler-2}, we can also get  \eqref{gen-r-abk}  by interchanging the order of summation in \eqref{0101-1} as follows.
\begin{align*}
\sum_{\pi\in\mathcal{R}_{a,b,k}}u^{\ell_a(\pi)}v^{\ell_b(\pi)}q^{|\pi|}&=\sum_{h\geq 0}\frac{v^hq^{k{h\choose 2}+bh}}{(q^k;q^k)_h}\sum_{m\geq h}\frac{u^{m-h}q^{(m-h)a}}{(q^k;q^k)_{m-h}}\\
&=\sum_{h\geq 0}\frac{v^hq^{k{h\choose 2}+bh}}{(q^k;q^k)_h}\sum_{m\geq 0}\frac{u^{m}q^{ma}}{(q^k;q^k)_{m}}\\
&=(-vq^{b};q^k)_\infty\cdot\frac{1}{(uq^a;q^k)_\infty}\\
&=\frac{(-vq^{b};q^k)_\infty}{(uq^a;q^k)_\infty}.
\end{align*}

The proof is complete.   \qed

\subsection{Proof of \eqref{gen-rr-abk}}

In this subsection, we will give a proof of \eqref{gen-rr-abk} with a similar argument in Section 4.1. For $m\geq 1$, let $\mathcal{BR}'_{a,b,k}(m)$ be the set of partitions $\lambda=(\lambda_1,\lambda_2,\ldots,\lambda_m)$ in $\mathcal{R}'_{a,b,k}$ with $m$ parts
 such that
\begin{itemize}
\item[(1)] $\lambda_m=a$ or $b$;

\item[(2)] for $1\leq i<m$, $\lambda_i-\lambda_{i+1}\leq k$ with strict inequality if $\lambda_{i+1}\equiv b\pmod{k}$.
\end{itemize}

Assume that $\lambda$ is a partition in $\mathcal{BR}'_{a,b,k}(m)$. For $1\leq i<m$, if $\lambda_{i+1}=kj+a$, then we have $kj+a<\lambda_{i}\leq k(j+1)+a$, and so  $\lambda_{i}=kj+b$ or $k(j+1)+a$. If $\lambda_{i+1}=kj+b$, then we have $kj+b\leq \lambda_{i}< k(j+1)+b$, and so  $\lambda_{i}=kj+b$ or $k(j+1)+a$. So, the number of partitions in $\mathcal{BR}'_{a,b,k}(m)$ is $2^m$. For example, the number of partitions in $\mathcal{BR}'_{a,b,k}(3)$ is $2^3=8$.
\[(b,b,a),(k+a,b,a),(k+b,k+a,a),(2k+a,k+a,a),\]
\[(b,b,b),(k+a,b,b),(k+b,k+a,b),(2k+a,k+a,b).\]

Set \[\mathcal{BR}'_{a,b,k}=\bigcup_{m\geq 1}\mathcal{BR}'_{a,b,k}(m).\]
 Obviously, $\mathcal{BR}'_{a,b,k}$ is the basis of $\mathcal{R}'_{a,b,k}$. So, we get
\begin{lem}
 $\mathcal{R}_{a,b,k}$ is a separable integer partition class with modulus $k$.
\end{lem}

For $m\geq 2$, let $\lambda$ be a partition in $\mathcal{BR}'_{a,b,k}(m)$. Clearly, we have
\[b\leq l(\lambda)\leq k(m-1)+a.\]

For $m\geq 2$ and $0\leq h\leq m-2$, let $\mathcal{BR}'_{a,b,k}(m,h,a)$ (resp. $\mathcal{BR}'_{a,b,k}(m,h,b)$) be the set of partitions in $\mathcal{BR}'_{a,b,k}(m)$ with the largest part $k(h+1)+a$ (resp. $kh+b$).

\begin{thm}\label{thm-gen-rr-abk-basis}
For $m\geq 2$ and $0\leq h\leq m-2$, we have
\begin{equation}\label{eqn-gen-rr-abk-basis-a}
\sum_{\lambda\in\mathcal{BR}'_{a,b,k}(m,h,a)}u^{\ell_a(\lambda)}v^{\ell_b(\lambda)}q^{|\lambda|}=(uq^a+vq^b)u^{h+1}v^{m-h-2}q^{(m-1)b+k{{h+1}\choose 2}+(k-b+a)(h+1)}{{m-2}\brack{h}}_k,
\end{equation}
and
\begin{equation*}\label{eqn-gen-rrr-abk-basis-b}
\sum_{\lambda\in\mathcal{BR}'_{a,b,k}(m,h,b)}u^{\ell_a(\lambda)}v^{\ell_b(\lambda)}q^{|\lambda|}=(uq^a+vq^b)u^{h}v^{m-h-1}q^{(m-1)b+k{{h+1}\choose 2}+(k-b+a)h}{{m-2}\brack{h}}_k.
\end{equation*}
\end{thm}

\pf Let $\lambda$ be a partition in $\mathcal{BR}'_{a,b,k}(m,h,a)$. By definition, we have $l(\lambda)=k(h+1)+a$. If we change the largest part $k(h+1)+a$ of $\lambda$ to $kh+b$, then we get a partition in $\mathcal{BR}_{a,b,k}(m,h,b)$, and vice versa. This implies that
\[\sum_{\lambda\in\mathcal{BR}'_{a,b,k}(m,h,a)}u^{\ell_a(\lambda)}v^{\ell_b(\lambda)}q^{|\lambda|}
=uv^{-1}q^{k-b+a}\sum_{\lambda\in\mathcal{BR}'_{a,b,k}(m,h,b)}u^{\ell_a(\lambda)}v^{\ell_b(\lambda)}q^{|\lambda|}.\]  We find that it suffices to show \eqref{eqn-gen-rr-abk-basis-a}.

For a partition $\lambda$ in $\mathcal{BR}_{a,b,k}(m,h,a)$, it is clear that the smallest part of $\lambda$ is $a$ or $b$ and there exist parts
\[k(h+1)+a,kh+a,\ldots,k+a\]
in $\lambda$. If we remove the smallest part and  the parts $k(h+1)+a,kh+a,\ldots,k+a$ from $\lambda$,
then we get a partition
with $m-h-2$ parts not exceeding ${kh+b}$ and equivalent to $b$ modulo $k$. The generating function for the partitions with $m-h-2$ parts not exceeding ${kh+b}$ and equivalent to $b$ modulo $k$ is
\[q^{b(m-h-2)}{{m-2}\brack{h}}_k.\]

Then, we get
\begin{align*}
\sum_{\lambda\in\mathcal{BR}'_{a,b,k}(m,h,a)}u^{\ell_a(\lambda)}v^{\ell_b(\lambda)}q^{|\lambda|}=(uq^a+vq^b)u^{h+1}q^{k{{h+2}\choose 2}+a(h+1)}\cdot v^{m-h-2}q^{b(m-h-2)}{{m-2}\brack{h}}_k,
\end{align*}
and so \eqref{eqn-gen-rr-abk-basis-a} is valid. This completes the proof.    \qed

Then, we give the generating function for the partitions in $\mathcal{BR}'_{a,b,k}(m)$.
\begin{thm}\label{rr-thm-000}
For $m\geq 1$,
\begin{equation}\label{new-002}
\sum_{\lambda\in\mathcal{BR}'_{a,b,k}(m)}u^{\ell_a(\lambda)}v^{\ell_b(\lambda)}q^{|\lambda|}=\sum_{h=0}^mu^hv^{m-h}q^{mb+k{h\choose 2}+(a-b)h}{{m}\brack{h}}_k.
\end{equation}
\end{thm}

\pf For $m=1$, it is clear that there are two partitions $(a)$ and $(b)$ in $\mathcal{BR}'_{a,b,k}(1)$, and so
\begin{equation*}
\sum_{\lambda\in\mathcal{BR}'_{a,b,k}(1)}u^{\ell_a(\lambda)}v^{\ell_b(\lambda)}q^{|\lambda|}=uq^a+vq^b,
\end{equation*}
which agrees with \eqref{new-002} for $m=1$.

For $m\geq 2$, in order to show \eqref{new-002}, it suffices to prove that
\begin{equation}\label{new-002-proof-1}
\sum_{\lambda\in\mathcal{BR}'_{a,b,k}(m)}u^{\ell_a(\lambda)}v^{\ell_b(\lambda)}q^{|\lambda|}=\sum_{h=0}^{m-1}(uq^a+vq^{b})u^hv^{m-h-1}q^{(m-1)b+k{h\choose 2}+(k-b+a)h}{{m-1}\brack{h}}_k,
\end{equation}
and
\begin{equation}\label{new-002-proof-2}
\sum_{h=0}^{m-1}(uq^a+vq^{b})u^hv^{m-h-1}q^{(m-1)b+k{h\choose 2}+(k-b+a)h}{{m-1}\brack{h}}_k=\sum_{h=0}^mu^hv^{m-h}q^{mb+k{h\choose 2}+(a-b)h}{{m}\brack{h}}_k.
\end{equation}

We first show \eqref{new-002-proof-1}. It follows from Theorem \ref{thm-gen-rr-abk-basis} that
\begin{equation}\label{use-111}
\begin{split}
&\quad\sum_{\lambda\in\mathcal{BR}'_{a,b,k}(m)}u^{\ell_a(\lambda)}v^{\ell_b(\lambda)}q^{|\lambda|}\\
&=\sum_{h=0}^{m-2}\left(\sum_{\lambda\in\mathcal{BR}'_{a,b,k}(m,h,a)}u^{\ell_a(\lambda)}v^{\ell_b(\lambda)}q^{|\lambda|}+\sum_{\lambda\in\mathcal{BR}'_{a,b,k}(m,h,b)}u^{\ell_a(\lambda)}v^{\ell_b(\lambda)}q^{|\lambda|}\right)\\
&=\sum_{h=0}^{m-2}\left((uq^a+vq^b)u^{h+1}v^{m-h-2}q^{(m-1)b+k{{h+1}\choose 2}+(k-b+a)(h+1)}{{m-2}\brack{h}}_k\right.\\
&\left.\qquad\qquad+(uq^a+vq^b)u^{h}v^{m-h-1}q^{(m-1)b+k{{h+1}\choose 2}+(k-b+a)h}{{m-2}\brack{h}}_k\right).
\end{split}
\end{equation}

Appealing to \eqref{bin-new-r-1}, we have
\begin{align*}
&\quad\sum_{h=0}^{m-2}\left(u^{h+1}v^{-h-1}q^{k{{h+1}\choose 2}+(k-b+a)(h+1)}{{m-2}\brack{h}}_k+u^{h}v^{-h}q^{k{{h+1}\choose 2}+(k-b+a)h}{{m-2}\brack{h}}_k\right)\\
&=\sum_{h=1}^{m-1}u^hv^{-h}q^{k{h\choose 2}+(k-b+a)h}{{m-2}\brack{h-1}}_k+\sum_{h=0}^{m-2}u^{h}v^{-h}q^{k{{h+1}\choose 2}+(k-b+a)h}{{m-2}\brack{h}}_k\\
&=u^{m-1}v^{-(m-1)}q^{k{{m-1}\choose 2}+(k-b+a)(m-1)}\\
&\qquad+\sum_{h=1}^{m-2}u^hv^{-h}q^{k{h\choose 2}+(k-b+a)h}\left({{m-2}\brack{h-1}}_k+q^{kh}{{m-2}\brack{h}}_k\right)+1\\
&=\sum_{h=0}^{m-1}u^hv^{-h}q^{k{h\choose 2}+(k-b+a)h}{{m-1}\brack{h}}_k.
\end{align*}
Combining with \eqref{use-111}, we arrive at \eqref{new-002-proof-1}.

It remains to show \eqref{new-002-proof-2}. Again by \eqref{bin-new-r-1}, we get
\begin{align*}
&\quad\sum_{h=0}^{m-1}(uq^a+vq^{b})u^hv^{m-h-1}q^{(m-1)b+k{h\choose 2}+(k-b+a)h}{{m-1}\brack{h}}_k\\
&=\sum_{h=0}^{m-1}u^{h+1}v^{m-h-1}q^{mb+k{h+1\choose 2}+(a-b)(h+1)}{{m-1}\brack{h}}_k+\sum_{h=0}^{m-1}u^hv^{m-h}q^{mb+k{h\choose 2}+kh+(a-b)h}{{m-1}\brack{h}}_k\\
&=\sum_{h=1}^{m}u^{h}v^{m-h}q^{mb+k{h\choose 2}+(a-b)h}{{m-1}\brack{h-1}}_k+\sum_{h=0}^{m-1}u^hv^{m-h}q^{mb+k{h\choose 2}+kh+(a-b)h}{{m-1}\brack{h}}_k\\
&=u^mq^{mb+k{m\choose 2}+(a-b)m}+\sum_{h=1}^{m-1}u^hv^{m-h}q^{mb+k{h\choose 2}+(a-b)h}\left({{m-1}\brack{h-1}}_k+q^{kh}{{m-1}\brack{h}}_k\right)+v^mq^{mb}\\
&=\sum_{h=0}^mu^hv^{m-h}q^{mb+k{h\choose 2}+(a-b)h}{{m}\brack{h}}_k.
\end{align*}
Hence, \eqref{new-002-proof-2} is verified.   This completes the proof.  \qed

Now, we are in a position to give a proof of \eqref{gen-rr-abk}.

{\noindent \bf Proof of \eqref{gen-rr-abk}.} In view of Theorem \ref{rr-thm-000}, we deduce that
\begin{align}
\sum_{\pi\in\mathcal{R}'_{a,b,k}}u^{\ell_a(\pi)}v^{\ell_b(\pi)}q^{|\pi|}&=1+\sum_{m\geq 1}\frac{1}{(q^k;q^k)_m}\sum_{\lambda\in\mathcal{BR}'_{a,b,k}(m)}u^{\ell_a(\lambda)}v^{\ell_b(\lambda)}q^{|\lambda|}\nonumber\\
&=\sum_{m\geq 0}\frac{1}{(q^k;q^k)_m}\sum_{h=0}^mu^hv^{m-h}q^{mb+k{h\choose 2}+(a-b)h}{{m}\brack{h}}_k.\label{rr-0101-1}
\end{align}

Summing up  the $h$-sum in \eqref{rr-0101-1} by letting $q\rightarrow q^{k}$ and $t=-uv^{-1}q^{a-b}$ in \eqref{h-sum}, we get
\begin{align*}
\sum_{\pi\in\mathcal{R}'_{a,b,k}}u^{\ell_a(\pi)}v^{\ell_b(\pi)}q^{|\pi|}&=\sum_{m\geq 0}\frac{v^mq^{mb}}{(q^k;q^k)_m}\sum_{h=0}^m(uv^{-1}q^{a-b})^hq^{k{h\choose 2}}{m\brack h}_k\\
&=\sum_{m\geq 0}\frac{v^mq^{mb}}{(q^k;q^k)_m}(-uv^{-1}q^{a-b};q^k)_m\\
&=\frac{(-uq^{a};q^k)_\infty}{(vq^b;q^k)_\infty},
\end{align*}
where the final equation follows from \eqref{chang-1}. We arrive at \eqref{gen-rr-abk}.

With the aid of \eqref{Euler-1} and \eqref{Euler-2}, we can also get  \eqref{gen-rr-abk}  by interchanging the order of summation in \eqref{rr-0101-1} as follows.
\begin{align*}
\sum_{\pi\in\mathcal{R}'_{a,b,k}}u^{\ell_a(\pi)}v^{\ell_b(\pi)}q^{|\pi|}
&=\sum_{h\geq 0}\frac{u^hq^{k{h\choose 2}+ah}}{(q^k;q^k)_h}\sum_{m\geq h}\frac{v^{m-h}q^{(m-h)b}}{(q^k;q^k)_{m-h}}\\
&=\sum_{h\geq 0}\frac{u^hq^{k{h\choose 2}+ah}}{(q^k;q^k)_h}\sum_{m\geq 0}\frac{v^{m}q^{mb}}{(q^k;q^k)_{m}}\\
&=(-uq^{a};q^k)_\infty\cdot\frac{1}{(vq^b;q^k)_\infty}\\
&=\frac{(-uq^{a};q^k)_\infty}{(vq^b;q^k)_\infty}.
\end{align*}

The proof is  complete.        \qed

\section{Proof of \eqref{gen-d-abk} and \eqref{gen-dd-abk}}

In this section, we will study $\mathcal{D}_{a,b,k}$ and $\mathcal{D}'_{a,b,k}$ in terms of separable integer partition classes with modulus $k$, and then give proofs of \eqref{gen-d-abk} and \eqref{gen-dd-abk} in Section 5.1 and 5.2 respectively.

\subsection{Proof of \eqref{gen-d-abk}}

For $m\geq 1$, let $\mathcal{BD}_{a,b,k}(m)$ be the set of partitions $\lambda=(\lambda_1,\lambda_2,\ldots,\lambda_m)$ in $\mathcal{D}_{a,b,k}$ with $m$ parts
 such that
 \begin{itemize}
\item[(1)] $\lambda_m=a$ or $b$;

\item[(2)] for $1\leq i<m$, $\lambda_i-\lambda_{i+1}\leq 2k$ with strict inequality if $\pi_{i+1}\equiv a\pmod{k}$.

\end{itemize}

Assume that $\lambda$ is a partition in $\mathcal{BD}_{a,b,k}(m)$. For $1\leq i<m$, if $\lambda_{i+1}=kj+a$, then we have $k(j+1)+a\leq\lambda_{i}<k(j+2)+a$, and so  $\lambda_{i}=k(j+1)+a$ or $k(j+1)+b$. If $\lambda_{i+1}=kj+b$, then we have $k(j+1)+b<\lambda_{i}\leq k(j+2)+b$, and so  $\lambda_{i}=k(j+2)+a$ or $k(j+2)+b$. So, the number of partitions in $\mathcal{BD}_{a,b,k}(m)$ is $2^m$. For example, the number of partitions in $\mathcal{BD}_{a,b,k}(3)$ is $2^3=8$.
\[(2k+a,k+a,a),(2k+b,k+a,a),(3k+a,k+b,a),(3k+b,k+b,a),\]
\[(3k+a,2k+a,b),(3k+b,2k+a,b),(4k+a,2k+b,b),(4k+b,2k+b,b).\]

\begin{lem}
 $\mathcal{D}_{a,b,k}$ is a separable integer partition class with modulus $k$.
\end{lem}

\pf  Set \[\mathcal{BD}_{a,b,k}=\bigcup_{m\geq 1}\mathcal{BD}_{a,b,k}(m).\]
 Obviously, $\mathcal{BD}_{a,b,k}$ is the basis of $\mathcal{D}_{a,b,k}$. This completes the proof.   \qed

For $m\geq 1$, let $\lambda$ be a partition in $\mathcal{BD}_{a,b,k}(m)$. Clearly, we have
we have
\[k(m-1)+a\leq l(\lambda) \leq 2k(m-1)+b.\]

For $m\geq 1$ and $0\leq h\leq m-1$, let $\mathcal{BD}_{a,b,k}(m,h,a)$ (resp. $\mathcal{BD}_{a,b,k}(m,h,b)$) be the set of partitions in $\mathcal{BD}_{a,b,k}(m)$ with the largest part $kh+k(m-1)+a$ (resp. $kh+k(m-1)+b$). The generating functions for the partitions in $\mathcal{BD}_{a,b,k}(m,h,a)$ and $\mathcal{BD}_{a,b,k}(m,h,b)$ respectively are given below.
\begin{thm}\label{thm-gen-d-abk-basis}
For $m\geq 1$ and $0\leq h\leq m-1$, we have
\begin{equation}\label{eqn-gen-d-abk-basis-a}
\sum_{\lambda\in\mathcal{BD}_{a,b,k}(m,h,a)}u^{\ell_a(\lambda)}v^{\ell_b(\lambda)}q^{|\lambda|}=u^{m-h}v^hq^{k{m\choose 2}+k{h+1\choose 2}+ma+(b-a)h}{{m-1}\brack{h}}_k,
\end{equation}
and
\begin{equation*}\label{eqn-gen-d-abk-basis-b}
\sum_{\lambda\in\mathcal{BD}_{a,b,k}(m,h,b)}u^{\ell_a(\lambda)}v^{\ell_b(\lambda)}q^{|\lambda|}
=u^{m-h-1}v^{h+1}q^{k{m\choose 2}+k{h+1\choose 2}+ma+(b-a)(h+1)}{{m-1}\brack{h}}_k.
\end{equation*}
\end{thm}

Before proving Theorem \ref{thm-gen-d-abk-basis}, we establish the following two relations.
\begin{lem}\label{lem-d-ab-0-1}
For $m\geq 1$ and $0\leq h\leq m-1$,
\begin{equation}\label{recur-d-ab-1}
\sum_{\lambda\in\mathcal{BD}_{a,b,k}(m,h,b)}u^{\ell_a(\lambda)}v^{\ell_b(\lambda)}q^{|\lambda|}
=u^{-1}vq^{b-a}\sum_{\lambda\in\mathcal{BD}_{a,b,k}(m,h,a)}u^{\ell_a(\lambda)}v^{\ell_b(\lambda)}q^{|\lambda|}.
\end{equation}
\end{lem}

\pf Let $\lambda$ be a partition in $\mathcal{BD}_{a,b,k}(m,h,b)$. By definition, we have $l(\lambda)=kh+k(m-1)+b$. If we change the largest part $kh+k(m-1)+b$ of $\lambda$ to $kh+k(m-1)+a$, then we get a partition in $\mathcal{BD}_{a,b,k}(m,h,a)$, and vice versa. This implies that \eqref{recur-d-ab-1} holds, and so the proof is complete.  \qed

\begin{lem}\label{lem-d-ab-recur-1}
For $m\geq 2$ and $1\leq h\leq m-1$,
\begin{equation}\label{recur-d-ab-2}
\begin{split}
&\quad\sum_{\lambda\in\mathcal{BD}_{a,b,k}(m+1,h,a)}u^{\ell_a(\lambda)}v^{\ell_b(\lambda)}q^{|\lambda|}\\
&=uq^{kh+km+a}\left(\sum_{\lambda\in\mathcal{BD}_{a,b,k}(m,h,a)}u^{\ell_a(\lambda)}v^{\ell_b(\lambda)}q^{|\lambda|}+\sum_{\lambda\in\mathcal{BD}_{a,b,k}(m,h-1,b)}u^{\ell_a(\lambda)}v^{\ell_b(\lambda)}q^{|\lambda|}\right).
\end{split}
\end{equation}
\end{lem}

\pf Let $\lambda$ be a partition in $\mathcal{BD}_{a,b,k}(m+1,h,a)$. By definition, we have $l(\lambda)=kh+km+a$. If we remove the largest part $kh+km+a$ from $\lambda$, then we get a  partition in $\mathcal{BD}_{a,b,k}(m)$ with the largest part
being $kh+k(m-1)+a$ or $k(h-1)+k(m-1)+b$, and vice versa. It yields that \eqref{recur-d-ab-2} is valid. This completes the proof.  \qed

Now, we proceed to show Theorem \ref{thm-gen-d-abk-basis}.

{\noindent \bf Proof of  Theorem \ref{thm-gen-d-abk-basis}.} It follows from Lemma \ref{lem-d-ab-0-1} that we just need to show \eqref{eqn-gen-d-abk-basis-a}. We will prove \eqref{eqn-gen-d-abk-basis-a} by induction on $m$.

For $m=1$, there is only one partition $(a)$ in $\mathcal{BD}_{a,b,k}(1,0,a)$, which yields
\[\sum_{\lambda\in\mathcal{BD}_{a,b,k}(1,0,a)}u^{\ell_a(\lambda)}v^{\ell_b(\lambda)}q^{|\lambda|}=uq^a.\]

 For $m=2$, it can be checked that there is only one partition $(k+a,a)$ in $\mathcal{BD}_{a,b,k}(2,0,a)$ and there is only one partition $(2k+a,b)$ in $\mathcal{BD}_{a,b,k}(2,1,a)$. This implies that
 \[\sum_{\lambda\in\mathcal{BD}_{a,b,k}(2,0,a)}u^{\ell_a(\lambda)}v^{\ell_b(\lambda)}q^{|\lambda|}=u^2q^{k+2a},\]
 and
 \[\sum_{\lambda\in\mathcal{BD}_{a,b,k}(2,1,a)}u^{\ell_a(\lambda)}v^{\ell_b(\lambda)}q^{|\lambda|}=uvq^{2k+a+b}.\]

For $m\geq 2$, assume that \eqref{eqn-gen-d-abk-basis-a} holds for $m$. We consider the following three cases.

Case 1: $h=0$. The only one partition in $\mathcal{BD}_{a,b,k}(m+1,0,a)$ is
\[(km+a,k(m-1)+a,\ldots,a).\]
 So, we have
\[\sum_{\lambda\in\mathcal{BD}_{a,b,k}(m+1,0,a)}u^{\ell_a(\lambda)}v^{\ell_b(\lambda)}q^{|\lambda|}=u^{m+1}q^{k{{m+1}\choose 2}+(m+1)a}.\]

Case 2: $h=m$. The only one partition in $\mathcal{BD}_{a,b,k}(m+1,m,a)$ is
\[(2km+a,2k(m-1)+b,\ldots,b).\]
 So, we have
\[\sum_{\lambda\in\mathcal{BD}_{a,b,k}(m+1,m,a)}u^{\ell_a(\lambda)}v^{\ell_b(\lambda)}q^{|\lambda|}=uv^{m}q^{2k{{m+1}\choose 2}+a+mb}.\]

Case 3: $1\leq h\leq m-1$. Under the assumption that \eqref{eqn-gen-d-abk-basis-a} holds for $m$, then by Lemma \ref{lem-d-ab-0-1} and Lemma \ref{lem-d-ab-recur-1}, we get
\begin{align*}
\quad&\sum_{\lambda\in\mathcal{BD}_{a,b,k}(m+1,h,a)}u^{\ell_a(\lambda)}v^{\ell_b(\lambda)}q^{|\lambda|}\\
&=uq^{kh+km+a}\left(\sum_{\lambda\in\mathcal{BD}_{a,b,k}(m,h,a)}u^{\ell_a(\lambda)}v^{\ell_b(\lambda)}q^{|\lambda|}
+\sum_{\lambda\in\mathcal{BD}_{a,b,k}(m,h-1,b)}u^{\ell_a(\lambda)}v^{\ell_b(\lambda)}q^{|\lambda|}\right)\\
&=uq^{kh+km+a}\left(\sum_{\lambda\in\mathcal{BD}_{a,b,k}(m,h,a)}u^{\ell_a(\lambda)}v^{\ell_b(\lambda)}q^{|\lambda|}
+u^{-1}vq^{b-a}\sum_{\lambda\in\mathcal{BD}_{a,b,k}(m,h-1,a)}u^{\ell_a(\lambda)}v^{\ell_b(\lambda)}q^{|\lambda|}\right)\\
&=uq^{kh+km+a}\left(u^{m-h}v^hq^{k{m\choose 2}+k{h+1\choose 2}+ma+(b-a)h}{{m-1}\brack{h}}_k\right.\\
&\qquad\qquad\qquad\qquad\left.+u^{m-h}v^hq^{k{m\choose 2}+k{h\choose 2}+ma+(b-a)h}{{m-1}\brack{h-1}}_k\right)\\
&=u^{m+1-h}v^hq^{k{m+1\choose 2}+k{h+1\choose 2}+(m+1)a+(b-a)h}\left(q^{kh}{{m-1}\brack{h}}_k+{{m-1}\brack{h-1}}_k\right)\\
&=u^{m+1-h}v^hq^{k{m+1\choose 2}+k{h+1\choose 2}+(m+1)a+(b-a)h}{{m}\brack{h}}_k,
\end{align*}
where the final equation follows from \eqref{bin-new-r-1}.

 We conclude that \eqref{eqn-gen-d-abk-basis-a} also holds for $m+1$. This completes the proof.  \qed

Then, we give the generating function for the partitions in $\mathcal{BD}_{a,b,k}(m)$.

\begin{thm}\label{thm-gen-abk-m-0}
For $m\geq 1$,
\begin{equation}\label{new-001}
\sum_{\lambda\in\mathcal{BD}_{a,b,k}(m)}u^{\ell_a(\lambda)}v^{\ell_b(\lambda)}q^{|\lambda|}
=\sum^m_{h=0}u^{m-h}v^hq^{k{m\choose 2}+k{h\choose 2}+ma+(b-a)h}{{m}\brack{h}}_k.
\end{equation}
\end{thm}

\pf  For $m=1$, there are two partitions
in $\mathcal{BD}_{a,b,k}(1)$, which are $(a)$ and $(b)$. It yields that
\begin{equation*}
\sum_{\lambda\in\mathcal{BD}_{a,b,k}(1)}u^{\ell_a(\lambda)}v^{\ell_b(\lambda)}q^{|\lambda|}=uq^a+vq^b,
\end{equation*}
which agrees with \eqref{new-001} for $m=1$.

For $m\geq 2$, it follows from Lemma \ref{lem-d-ab-recur-1} that
\begin{align*}
&\quad \sum_{\lambda\in\mathcal{BD}_{a,b,k}(m)}u^{\ell_a(\lambda)}v^{\ell_b(\lambda)}q^{|\lambda|}\\
&=\sum^{m-1}_{h=0}\left(\sum_{\lambda\in\mathcal{BD}_{a,b,k}(m,h,a)}u^{\ell_a(\lambda)}v^{\ell_b(\lambda)}q^{|\lambda|}
+\sum_{\lambda\in\mathcal{BD}_{a,b,k}(m,h,b)}u^{\ell_a(\lambda)}v^{\ell_b(\lambda)}q^{|\lambda|}\right)\\
&=\sum^{m-1}_{h=0}\sum_{\lambda\in\mathcal{BD}_{a,b,k}(m,h,a)}u^{\ell_a(\lambda)}v^{\ell_b(\lambda)}q^{|\lambda|}
+\sum^{m}_{h=1}\sum_{\lambda\in\mathcal{BD}_{a,b,k}(m,h-1,b)}u^{\ell_a(\lambda)}v^{\ell_b(\lambda)}q^{|\lambda|}\\
&=\sum_{\lambda\in\mathcal{BD}_{a,b,k}(m,0,a)}u^{\ell_a(\lambda)}v^{\ell_b(\lambda)}q^{|\lambda|}\\
&\qquad+\sum^{m-1}_{h=1}\left(\sum_{\lambda\in\mathcal{BD}_{a,b,k}(m,h,a)}u^{\ell_a(\lambda)}v^{\ell_b(\lambda)}q^{|\lambda|}
+\sum_{\lambda\in\mathcal{BD}_{a,b,k}(m,h-1,b)}u^{\ell_a(\lambda)}v^{\ell_b(\lambda)}q^{|\lambda|}\right)\\
&\qquad+\sum_{\lambda\in\mathcal{BD}_{a,b,k}(m,m-1,b)}u^{\ell_a(\lambda)}v^{\ell_b(\lambda)}q^{|\lambda|}\\
&=\sum_{\lambda\in\mathcal{BD}_{a,b,k}(m,0,a)}u^{\ell_a(\lambda)}v^{\ell_b(\lambda)}q^{|\lambda|}\\
&\qquad+\sum^{m-1}_{h=1}(uq^{kh+km+a})^{-1}\sum_{\lambda\in\mathcal{BD}_{a,b,k}(m+1,h,a)}u^{\ell_a(\lambda)}v^{\ell_b(\lambda)}q^{|\lambda|}+\sum_{\lambda\in\mathcal{BD}_{a,b,k}(m,m-1,b)}u^{\ell_a(\lambda)}v^{\ell_b(\lambda)}q^{|\lambda|}.
\end{align*}
Combining with Theorem \ref{thm-gen-d-abk-basis}, we get
\begin{align*}
&\quad \sum_{\lambda\in\mathcal{BD}_{a,b,k}(m)}u^{\ell_a(\lambda)}v^{\ell_b(\lambda)}q^{|\lambda|}\\
&=u^mq^{k{m\choose 2}+ma}+\sum^{m-1}_{h=1}u^{m-h}v^hq^{k{m\choose 2}+k{h\choose 2}+ma+(b-a)h}{{m}\brack{h}}_k+v^mq^{2k{m\choose 2}+mb}\\
&=\sum^{m}_{h=0}u^{m-h}v^hq^{k{m\choose 2}+k{h\choose 2}+ma+(b-a)h}{{m}\brack{h}}_k.
\end{align*}

We arrive at \eqref{new-001}, and thus the proof is complete.  \qed

Now, we are in a position to give the proof of \eqref{gen-d-abk}.

{\noindent \bf Proof of \eqref{gen-d-abk}.} Appealing to Theorem \ref{thm-gen-abk-m-0}, we get
\begin{align}
\sum_{\pi\in\mathcal{D}_{a,b,k}}u^{\ell_a(\pi)}v^{\ell_b(\pi)}q^{|\pi|}&=1+\sum_{m\geq 1}\frac{1}{(q^k;q^k)_m}\sum_{\lambda\in\mathcal{BD}_{a,b,k}(m)}u^{\ell_a(\lambda)}v^{\ell_b(\lambda)}q^{|\lambda|}\nonumber\\
&=\sum_{m\geq 0}\frac{1}{(q^k;q^k)_m}\sum^{m}_{h=0}u^{m-h}v^hq^{k{m\choose 2}+k{h\choose 2}+ma+(b-a)h}{{m}\brack{h}}_k.\label{abk-inter-double-sum}
\end{align}

Summing up  the $h$-sum in \eqref{abk-inter-double-sum} by letting $q\rightarrow q^{k}$ and $t=-u^{-1}vq^{b-a}$ in \eqref{h-sum}, we get
\begin{align}
\sum_{\pi\in\mathcal{D}_{a,b,k}}u^{\ell_a(\pi)}v^{\ell_b(\pi)}q^{|\pi|}&=\sum_{m\geq 0}\frac{u^mq^{k{m\choose 2}+ma}}{(q^k;q^k)_m}\sum_{h=0}^m(u^{-1}vq^{b-a})^hq^{k{h\choose 2}}{m\brack h}_k\nonumber\\
&=\sum_{m\geq 0}\frac{u^mq^{k{m\choose 2}+ma}}{(q^k;q^k)_m}(-u^{-1}vq^{b-a};q^k)_m.\nonumber
\end{align}

On the other hand, by interchanging the order of summation in \eqref{abk-inter-double-sum}, we can obtain that
\begin{align}
\sum_{\pi\in\mathcal{D}_{a,b,k}}u^{\ell_a(\pi)}v^{\ell_b(\pi)}q^{|\pi|}
&=\sum_{h\geq 0}\frac{v^hq^{k{h\choose 2}+(b-a)h}}{(q^k;q^k)_h}\sum_{m\geq h}\frac{u^{m-h}q^{k{m\choose 2}+ma}}{(q^k;q^k)_{m-h}}\nonumber\\
&=\sum_{h\geq 0}\frac{v^hq^{k{h\choose 2}+(b-a)h}}{(q^k;q^k)_h}\sum_{m\geq 0}\frac{u^mq^{k{{m+h}\choose 2}+(m+h)a}}{(q^k;q^k)_m}\nonumber\\
&=\sum_{h\geq 0}\frac{v^hq^{2k{h\choose 2}+bh}}{(q^k;q^k)_h}\sum_{m\geq 0}\frac{(uq^{kh+a})^mq^{k{m\choose 2}}}{(q^k;q^k)_m}\nonumber\\
&=\sum_{h\geq 0}\frac{v^hq^{2k{h\choose 2}+bh}}{(q^k;q^k)_h}(-uq^{kh+a};q^k)_\infty,\nonumber
\end{align}
where the final equation follows from \eqref{Euler-2} with  $q\rightarrow q^{k}$ and $t=uq^{kh+a}$. This completes the proof.  \qed

\subsection{Proof of \eqref{gen-dd-abk}}

In this subsection, we will give a proof of \eqref{gen-dd-abk} with a similar argument in Section 5.1.
For $m\geq 1$, let $\mathcal{BD}'_{a,b,k}(m)$ be the set of partitions $\lambda=(\lambda_1,\lambda_2,\ldots,\lambda_m)$ in $\mathcal{D}_{a,b,k}$ with $m$ parts
 such that
 \begin{itemize}
\item[(1)] $\lambda_m=a$ or $b$;

\item[(2)] for $1\leq i<m$, $\lambda_i-\lambda_{i+1}\leq 2k$ with strict inequality if $\pi_{i+1}\equiv b\pmod{k}$.

\end{itemize}

 Assume that $\lambda$ is a partition in $\mathcal{BD}'_{a,b,k}(m)$. For $1\leq i<m$, if $\lambda_{i+1}=kj+a$, then we have $k(j+1)+a<\lambda_{i}\leq k(j+2)+a$, and so  $\lambda_{i}=k(j+1)+b$ or $k(j+2)+a$. If $\lambda_{i+1}=kj+b$, then we have $k(j+1)+b\leq \lambda_{i}<k(j+2)+b$, and so  $\lambda_{i}=k(j+1)+b$ or $k(j+2)+a$. So, the number of partitions in $\mathcal{BD}'_{a,b,k}(m)$ is $2^m$. For example, the number of partitions in $\mathcal{BD}'_{a,b,k}(3)$ is $2^3=8$.
\[(2k+b,k+b,a),(3k+a,k+b,a),(3k+b,2k+a,a),(4k+a,2k+a,a),\]
\[(2k+b,k+b,b),(3k+a,k+b,b),(3k+b,2k+a,b),(4k+a,2k+a,b).\]

Set \[\mathcal{BD}'_{a,b,k}=\bigcup_{m\geq 1}\mathcal{BD}'_{a,b,k}(m).\]
 Obviously, $\mathcal{BD}'_{a,b,k}$ is the basis of $\mathcal{D}'_{a,b,k}$. So, we get
\begin{lem}
 $\mathcal{D}'_{a,b,k}$ is a separable integer partition class with modulus $k$.
\end{lem}

For $m\geq 2$, let $\lambda$ be a partition in $\mathcal{BD}'_{a,b,k}(m)$. Clearly, we have
\[k(m-1)+b\leq l(\lambda) \leq 2k(m-1)+a\]

For $m\geq 2$ and $0\leq h\leq m-2$, let $\mathcal{BD}'_{a,b,k}(m,h,a)$ (resp. $\mathcal{BD}'_{a,b,k}(m,h,b)$) be the set of partitions in $\mathcal{BD}'_{a,b,k}(m)$ with the largest part $k(h+1)+k(m-1)+a$ (resp. $kh+k(m-1)+b$). Similarly to Lemma \ref{lem-d-ab-0-1} and Lemma \ref{lem-d-ab-recur-1}, we have
\begin{lem}\label{lem-dd-ab-0-1} For $m\geq 2$ and $0\leq h\leq m-2$,
\[\sum_{\lambda\in\mathcal{BD}'_{a,b,k}(m,h,a)}u^{\ell_a(\lambda)}v^{\ell_b(\lambda)}q^{|\lambda|}
=uv^{-1}q^{k-b+a}\sum_{\lambda\in\mathcal{BD}'_{a,b,k}(m,h,b)}u^{\ell_a(\lambda)}v^{\ell_b(\lambda)}q^{|\lambda|}.\]
\end{lem}
\begin{lem}\label{lem-dd-ab-recur-1}
For $m\geq 3$ and $1\leq h\leq m-2$,
\begin{align*}
&\quad \sum_{\lambda\in\mathcal{BD}'_{a,b,k}(m+1,h,b)}u^{\ell_a(\lambda)}v^{\ell_b(\lambda)}q^{|\lambda|}\\
&=vq^{kh+km+b}\left(\sum_{\lambda\in\mathcal{BD}'_{a,b,k}(m,h-1,a)}u^{\ell_a(\lambda)}v^{\ell_b(\lambda)}q^{|\lambda|}
+\sum_{\lambda\in\mathcal{BD}'_{a,b,k}(m,h,b)}u^{\ell_a(\lambda)}v^{\ell_b(\lambda)}q^{|\lambda|}\right).
\end{align*}
\end{lem}

Then, we give the generating functions for the partitions in $\mathcal{BD}_{a,b,k}(m,h,a)$ and $\mathcal{BD}_{a,b,k}(m,h,b)$ respectively.
\begin{thm}\label{thm-gen-dd-abk-basis}
For $m\geq 2$ and $0\leq h\leq m-2$, we have
\begin{align*}\label{eqn-gen-dd-abk-basis-b}
&\quad\sum_{\lambda\in\mathcal{BD}'_{a,b,k}(m,h,a)}u^{\ell_a(\lambda)}v^{\ell_b(\lambda)}q^{|\lambda|}\\
&=(uq^a+vq^b)u^{h+1}v^{m-h-2}q^{k{m\choose 2}+k{h+1\choose 2}+(m-1)b+(k-b+a)(h+1)}{{m-2}\brack{h}}_k,
\end{align*}
and
\begin{equation}\label{eqn-gen-dd-abk-basis-a}
\sum_{\lambda\in\mathcal{BD}'_{a,b,k}(m,h,b)}u^{\ell_a(\lambda)}v^{\ell_b(\lambda)}q^{|\lambda|}=(uq^a+vq^b)u^hv^{m-h-1}q^{k{m\choose 2}+k{h+1\choose 2}+(m-1)b+(k-b+a)h}{{m-2}\brack{h}}_k.
\end{equation}
\end{thm}
\pf It is clear from Lemma \ref{lem-dd-ab-0-1} that we just need to show \eqref{eqn-gen-dd-abk-basis-a}. We will prove \eqref{eqn-gen-dd-abk-basis-a} by induction on $m$.

For $m=2$, there are two partitions  in $\mathcal{BD}'_{a,b,k}(2,0,b)$, which are $(k+b,a)$ and $(k+b,b)$. It yields that
\[\sum_{\lambda\in\mathcal{BD}'_{a,b,k}(2,0,b)}u^{\ell_a(\lambda)}v^{\ell_b(\lambda)}q^{|\lambda|}=(uq^a+vq^b)vq^{k+b}.\]

For $m=3$, there are two partitions $(2k+b,k+b,a)$ and $(2k+b,k+b,b)$ in $\mathcal{BD}'_{a,b,k}(3,0,b)$ and there are two partitions $(3k+b,2k+a,a)$ and $(3k+b,2k+a,b)$ in $\mathcal{BD}'_{a,b,k}(3,1,b)$. This implies that
\[\sum_{\lambda\in\mathcal{BD}'_{a,b,k}(3,0,b)}u^{\ell_a(\lambda)}v^{\ell_b(\lambda)}q^{|\lambda|}=(uq^a+vq^b)v^2q^{3k+2b},\]
and
\[\sum_{\lambda\in\mathcal{BD}'_{a,b,k}(3,1,b)}u^{\ell_a(\lambda)}v^{\ell_b(\lambda)}q^{|\lambda|}=(uq^a+vq^b)v^2q^{5k+b+a}.\]

For $m\geq 3$, assume that \eqref{eqn-gen-dd-abk-basis-a} holds for $m$. We consider the following three cases.

Case 1: $h=0$. There are two partitions in $\mathcal{BD}'_{a,b,k}(m+1,0,b)$, which are
\[(km+b,k(m-1)+b,\ldots,k+b,a)\text{ and }(km+b,k(m-1)+b,\ldots,k+b,b).\]
 This implies that
\[\sum_{\lambda\in\mathcal{BD}'_{a,b,k}(m+1,0,b)}u^{\ell_a(\lambda)}v^{\ell_b(\lambda)}q^{|\lambda|}=(uq^a+vq^b)v^{m}q^{k{{m+1}\choose 2}+mb}.\]

Case 2: $h=m-1$. There are two partitions in $\mathcal{BD}'_{a,b,k}(m+1,m-1,b)$, which are
\[(k(2m-1)+b,k(2m-2)+a,\ldots,2k+a,a)\text{ and }(k(2m-1)+b,k(2m-2)+a,\ldots,2k+a,b).\]
 It yields that
\[\sum_{\lambda\in\mathcal{BD}'_{a,b,k}(m+1,m-1,b)}u^{\ell_a(\lambda)}v^{\ell_b(\lambda)}q^{|\lambda|}=(uq^a+vq^b)u^{m-1}vq^{k(m^2+m-1)+b+(m-1)a}.\]

Case 3: $1\leq h\leq m-2$. Under the assumption that \eqref{eqn-gen-dd-abk-basis-a} holds for $m$, then by Lemma \ref{lem-dd-ab-0-1} and Lemma \ref{lem-dd-ab-recur-1}, we get
\begin{align*}
\quad&\sum_{\lambda\in\mathcal{BD}'_{a,b,k}(m+1,h,b)}u^{\ell_a(\lambda)}v^{\ell_b(\lambda)}q^{|\lambda|}\\
&=vq^{kh+km+b}\left(\sum_{\lambda\in\mathcal{BD}'_{a,b,k}(m,h-1,a)}u^{\ell_a(\lambda)}v^{\ell_b(\lambda)}q^{|\lambda|}
+\sum_{\lambda\in\mathcal{BD}'_{a,b,k}(m,h,b)}u^{\ell_a(\lambda)}v^{\ell_b(\lambda)}q^{|\lambda|}\right)\\
&=vq^{kh+km+b}\left(uv^{-1}q^{k-b+a}\sum_{\lambda\in\mathcal{BD}'_{a,b,k}(m,h-1,b)}u^{\ell_a(\lambda)}v^{\ell_b(\lambda)}q^{|\lambda|}
+\sum_{\lambda\in\mathcal{BD}'_{a,b,k}(m,h,b)}u^{\ell_a(\lambda)}v^{\ell_b(\lambda)}q^{|\lambda|}\right)\\
&=vq^{kh+km+b}\left((uq^a+vq^b)u^hv^{m-h-1}q^{k{m\choose 2}+k{h\choose 2}+(m-1)b+(k-b+a)h}{{m-2}\brack{h-1}}_k\right.\\
&\qquad\qquad\qquad\qquad\left.+(uq^a+vq^b)u^hv^{m-h-1}q^{k{m\choose 2}+k{h+1\choose 2}+(m-1)b+(k-b+a)h}{{m-2}\brack{h}}_k\right)\\
&=(uq^a+vq^b)u^hv^{m-h}q^{k{m+1\choose 2}+k{h+1\choose 2}+mb+(k-b+a)h}\left({{m-2}\brack{h-1}}_k+q^{kh}{{m-2}\brack{h}}_k\right)\\
&=(uq^a+vq^b)u^hv^{m-h}q^{k{m+1\choose 2}+k{h+1\choose 2}+mb+(k-b+a)h}{{m-1}\brack{h}}_k,\\
\end{align*}
where the final equation follows from \eqref{bin-new-r-1}.

 We conclude that \eqref{eqn-gen-dd-abk-basis-a} is satisfied for $m+1$. The proof is complete.  \qed

Now, we proceed to give the generating function for the partitions in $\mathcal{BD}'_{a,b,k}(m)$.

\begin{thm}\label{thm-gen-abk-m-1}
For $m\geq 1$,
\begin{equation*}
\sum_{\lambda\in\mathcal{BR}'_{a,b,k}(m)}u^{\ell_a(\lambda)}v^{\ell_b(\lambda)}q^{|\lambda|}=\sum_{h=0}^mu^hv^{m-h}q^{k{m\choose 2}+k{h\choose 2}+mb+(a-b)h}{{m}\brack{h}}_k.
\end{equation*}
\end{thm}

\pf In view of \eqref{new-002-proof-2}, we find that it suffices to prove that for $m\geq 1$,
 \begin{equation}\label{new-rr-002}
\sum_{\lambda\in\mathcal{BD}'_{a,b,k}(m)}u^{\ell_a(\lambda)}v^{\ell_b(\lambda)}q^{|\lambda|}
=\sum^{m-1}_{h=0}(uq^a+vq^b)u^hv^{m-h-1}q^{k{m\choose 2}+k{h\choose 2}+(m-1)b+(k-b+a)h}{{m-1}\brack{h}}_k.
\end{equation}

For $m=1$, it is clear that
\begin{equation*}
\sum_{\lambda\in\mathcal{BD}'_{a,b,k}(1)}u^{\ell_a(\lambda)}v^{\ell_b(\lambda)}q^{|\lambda|}=uq^a+vq^b,
\end{equation*}
since there are two partitions $(a)$ and $(b)$ in $\mathcal{BD}'_{a,b,k}(1)$. This implies that \eqref{new-rr-002} holds for $m=1$.

 For $m=2$, by Theorem \ref{thm-gen-dd-abk-basis}, we get
 \begin{align*}
 \sum_{\lambda\in\mathcal{BD}'_{a,b,k}(2)}u^{\ell_a(\lambda)}v^{\ell_b(\lambda)}q^{|\lambda|}&=\sum_{\lambda\in\mathcal{BD}'_{a,b,k}(2,0,b)}u^{\ell_a(\lambda)}v^{\ell_b(\lambda)}q^{|\lambda|}+\sum_{\lambda\in\mathcal{BD}'_{a,b,k}(2,0,a)}u^{\ell_a(\lambda)}v^{\ell_b(\lambda)}q^{|\lambda|}\\
 &=(uq^a+vq^b)vq^{k+b}+(uq^a+vq^b)uq^{2k+a},
 \end{align*}
which agrees with \eqref{new-rr-002} for $m=2$.

 For $m\geq 3$, it follows from Lemma \ref{lem-dd-ab-recur-1} that
 \begin{align*}
&\quad \sum_{\lambda\in\mathcal{BD}'_{a,b,k}(m)}u^{\ell_a(\lambda)}v^{\ell_b(\lambda)}q^{|\lambda|}\\
&=\sum^{m-2}_{h=0}\left(\sum_{\lambda\in\mathcal{BD}'_{a,b,k}(m,h,a)}u^{\ell_a(\lambda)}v^{\ell_b(\lambda)}q^{|\lambda|}
+\sum_{\lambda\in\mathcal{BD}'_{a,b,k}(m,h,b)}u^{\ell_a(\lambda)}v^{\ell_b(\lambda)}q^{|\lambda|}\right)\\
&=\sum^{m-1}_{h=1}\sum_{\lambda\in\mathcal{BD}'_{a,b,k}(m,h-1,a)}u^{\ell_a(\lambda)}v^{\ell_b(\lambda)}q^{|\lambda|}
+\sum^{m-2}_{h=0}\sum_{\lambda\in\mathcal{BD}'_{a,b,k}(m,h,b)}u^{\ell_a(\lambda)}v^{\ell_b(\lambda)}q^{|\lambda|}\\
&=\sum_{\lambda\in\mathcal{BD}'_{a,b,k}(m,0,b)}u^{\ell_a(\lambda)}v^{\ell_b(\lambda)}q^{|\lambda|}\\
&\qquad +\sum^{m-2}_{h=1}\left(\sum_{\lambda\in\mathcal{BD}'_{a,b,k}(m,h-1,a)}u^{\ell_a(\lambda)}v^{\ell_b(\lambda)}q^{|\lambda|}
+\sum_{\lambda\in\mathcal{BD}'_{a,b,k}(m,h,b)}u^{\ell_a(\lambda)}v^{\ell_b(\lambda)}q^{|\lambda|}\right)\\
&\qquad+\sum_{\lambda\in\mathcal{BD}'_{a,b,k}(m,m-2,a)}u^{\ell_a(\lambda)}v^{\ell_b(\lambda)}q^{|\lambda|}\\
&=\sum_{\lambda\in\mathcal{BD}'_{a,b,k}(m,0,b)}u^{\ell_a(\lambda)}v^{\ell_b(\lambda)}q^{|\lambda|}\\
&\qquad+\sum^{m-2}_{h=1}(vq^{kh+km+b})^{-1}\sum_{\lambda\in\mathcal{BD}'_{a,b,k}(m+1,h,b)}u^{\ell_a(\lambda)}v^{\ell_b(\lambda)}q^{|\lambda|}
+\sum_{\lambda\in\mathcal{BD}'_{a,b,k}(m,m-2,a)}u^{\ell_a(\lambda)}v^{\ell_b(\lambda)}q^{|\lambda|}.
\end{align*}
 Combining with Theorem \ref{thm-gen-dd-abk-basis}, we get
\begin{align*}
&\quad\sum_{\lambda\in\mathcal{BD}'_{a,b,k}(m)}u^{\ell_a(\lambda)}v^{\ell_b(\lambda)}q^{|\lambda|}\\
&=(uq^a+vq^b)v^{m-1}q^{k{m\choose 2}+(m-1)b}\\
&\quad +\sum^{m-2}_{h=1}u^hv^{m-h-1}(uq^a+vq^b)q^{k{m\choose 2}+k{h\choose 2}+(m-1)b+(k-b+a)h}{{m-1}\brack{h}}_k\\
&\quad +(uq^a+vq^b)u^{m-1}q^{k{m\choose 2}+k{m-1\choose 2}+(m-1)b+(k-b+a)(m-1)}\\
&=\sum^{m-1}_{h=0}u^hv^{m-h-1}(uq^a+vq^b)q^{k{m\choose 2}+k{h\choose 2}+(m-1)b+(k-b+a)h}{{m-1}\brack{h}}_k.\\
\end{align*}
Hence, \eqref{new-rr-002} is valid for $m\geq 3$. This completes the proof.  \qed

Finally, we conclude this section with a proof of \eqref{gen-dd-abk}.

{\noindent \bf Proof of \eqref{gen-dd-abk}.}  Utilizing Theorem \ref{thm-gen-abk-m-1}, we get
\begin{align}
\sum_{\pi\in\mathcal{D}'_{a,b,k}}u^{\ell_a(\pi)}v^{\ell_b(\pi)}q^{|\pi|}&=1+\sum_{m\geq 1}\frac{1}{(q^k;q^k)_m}\sum_{\lambda\in\mathcal{BD}'_{a,b,k}(m)}u^{\ell_a(\lambda)}v^{\ell_b(\lambda)}q^{|\lambda|}\nonumber\\
&=\sum_{m\geq 0}\frac{1}{(q^k;q^k)_m}\sum_{h=0}^mu^hv^{m-h}q^{k{m\choose 2}+k{h\choose 2}+mb+(a-b)h}{{m}\brack{h}}_k.\label{abk-inter-double-sum-dd-1}
\end{align}

Summing up  the $h$-sum in \eqref{abk-inter-double-sum-dd-1} by letting $q\rightarrow q^{k}$ and $t=-uv^{-1}q^{a-b}$ in \eqref{h-sum}, we get
\begin{align}
\sum_{\pi\in\mathcal{D}_{a,b,k}}u^{\ell_a(\pi)}v^{\ell_b(\pi)}q^{|\pi|}
&=\sum_{m\geq 0}\frac{v^mq^{k{m\choose 2}+mb}}{(q^k;q^k)_m}\sum_{h=0}^m(uv^{-1}q^{a-b})^hq^{k{h\choose 2}}{m\brack h}_k\nonumber\\
&=\sum_{m\geq 0}\frac{v^mq^{k{m\choose 2}+mb}}{(q^k;q^k)_m}(-uv^{-1}q^{a-b};q^k)_m.\nonumber
\end{align}

On the other hand, by interchanging the order of summation in \eqref{abk-inter-double-sum-dd-1}, we can obtain that
\begin{align*}
\sum_{\pi\in\mathcal{D}'_{a,b,k}}u^{\ell_a(\pi)}v^{\ell_b(\pi)}q^{|\pi|}&=\sum_{m\geq 0}\frac{v^mq^{k{m\choose 2}+mb}}{(q^k;q^k)_m}\sum_{h=0}^mu^hv^{-h}q^{k{h\choose 2}+(a-b)h}{m\brack h}_k\\
&=\sum_{h\geq 0}\frac{u^hq^{k{h\choose 2}+ah}}{(q^k;q^k)_h}\sum_{m\geq h}\frac{v^{m-h}q^{k{m\choose 2}+(m-h)b}}{(q^k;q^k)_{m-h}}\\
&=\sum_{h\geq 0}\frac{u^hq^{k{h\choose 2}+ah}}{(q^k;q^k)_h}\sum_{m\geq o}\frac{v^mq^{k{{m+h}\choose 2}+mb}}{(q^k;q^k)_m}\\
&=\sum_{h\geq 0}\frac{u^hq^{2k{h\choose 2}+ah}}{(q^k;q^k)_h}\sum_{m\geq 0}\frac{(vq^{kh+b})^mq^{k{m\choose 2}}}{(q^k;q^k)_m}\\
&=\sum_{h\geq 0}\frac{u^hq^{2k{h\choose 2}+ah}}{(q^k;q^k)_h}(-vq^{kh+b};q^k)_\infty,
\end{align*}
 where the final equation follows from \eqref{Euler-2} with  $q\rightarrow q^{k}$ and $t=vq^{kh+b}$. This completes the proof.  \qed

{\noindent \bf \large Declarations}

{\noindent \bf Conflict of interest:} The authors have no competing interests to declare that are relevant to the content of this article.

\end{document}